\documentclass[10pt,openbib]{article}

\usepackage[T1]{fontenc}
\usepackage{amsfonts,amsthm}
\usepackage{amsmath}
\usepackage{amssymb}

\newcommand{\R}{\mathbb{R}}

\newcommand{\Z}{\mathbb{Z}}

\newcommand{\Gl}{\mathop{\mathrm{Gl}}\nolimits}

\newcommand{\comp}{\, \raisebox{2pt}{$\scriptstyle\circ \, $}}

\newcommand{\vvee}{\mbox{\tiny $\vee$}}
\allowdisplaybreaks

\begin{document}

\title{\textbf{Vector fields invariant under a linear}
\\ \textbf{action of a compact Lie group} } 
\author{Richard Cushman\footnotemark}
\addtocounter{footnote}{1}
\footnotetext{email: r.h.cushman@gmail.com 
}
\date{}
\maketitle

The main result of this note is to give a new proof of: the $C^{\infty}({\R }^n)^G$ module of smooth vector fields on ${\R }^n$, which are invariant under a linear action of a compact subgroup $G$ of $\Gl (n, \R)$ is generated by polynomial vector fields on ${\R }^n$ that are invariant under this action. This result is standard and its proof can be found in \cite{golubitsky-stewart-schaeffer} and \cite{po\`{e}naru}. \medskip 

Let 
\begin{equation}
\Phi : G \times {\R }^n \rightarrow {\R }^n : (g,x) \mapsto {\Phi }_g(x) = g \cdot x
\label{eq-s1one}
\end{equation}
be a linear action of a compact subgroup $G$ of $\Gl (n, \R )$. In other words, 
for every $g \in G$ the mapping ${\Phi }_g: {\R }^n \rightarrow {\R }^n: 
x \mapsto {\Phi }_g(x)$ is an element of $\Gl (n, \R )$. So 
$G \rightarrow \Gl (n, \R): g \mapsto {\Phi }_g$ is an $n$ dimensional 
representation of $G$. \medskip 

Our main result is \medskip 

\noindent \textbf{Theorem 1} The $C^{\infty}({\R }^n)^G$ module 
$\mathfrak{X}({\R }^n)^G$ of $G$ invariant smooth vector fields on 
${\R }^n$ is finitely generated by polynomial vector fields that generate 
the $\mathrm{P}({\R }^n, \R)^G$ module $\mathcal{P}\mathfrak{X}({\R }^n)$ of 
$G$ invariant polynomial vector fields on ${\R }^n$. \medskip 

We need some preparatory results. \medskip

The linear $G$ action $\Phi $ induces a linear $G$ action on the 
space $\mathrm{P}({\R }^n, \R)$ of real valued polynomials on ${\R }^n$ defined by \begin{displaymath}
{\Phi }^{\dagger} : G \times \mathrm{P}({\R }^n, \R ) \rightarrow 
\mathrm{P}({\R }^n, \R ): (g, p) \mapsto {\Phi }^{\dagger}_gp, 
\end{displaymath}
where $({\Phi }^{\dagger }_gp)(x) = p({\Phi }_{g^{-1}}(x))$ for every $(g,x) \in 
G \times {\R }^n$. The mapping ${\Phi }^{\dagger }$ is a linear action, because 
for every $g \in G$ the mapping $ p \mapsto {\Phi }^{\dagger}_gp$ is 
linear and for every $(h,x) \in G \times {\R }^n$
\begin{align*}
({\Phi }^{\dagger}_{gh}p)(x) & = p({\Phi}_{(gh)^{-1}}(x)) =
p( {\Phi }_{h^{-1}} ({\Phi }_{g^{-1}}(x)))  \\
& = ({\Phi }^{\dagger}_hp)({\Phi }_{g^{-1}}(x)) 
= {\Phi }^{\dagger}_g({\Phi }^{\dagger}_hp)(x). 
\end{align*} 

\noindent \textbf{Fact 2} The algebra $\mathrm{P}({\R }^n, \R )^G$ of 
$G$ invariant polynomials on ${\R}^n$ is finitely generated. \medskip 

\noindent \textbf{Proof.} See invariant theory in \cite{wikipedia}. \hfill $\square $ \medskip 

\noindent \textbf{Fact 3} The algebra $C^{\infty}({\R }^n)^G$ of smooth 
$G$ invariant functions on ${\R }^n$ is equal to $C^{\infty}({\R }^k)$, 
where ${\R }^k$ has coordinates $p_1, \ldots , p_k$, which are the 
generators of $\mathrm{P}({\R }^n, \R )^G$. \medskip 

\noindent \textbf{Proof.} See Schwarz \cite{schwarz}. \hfill $\square $ \medskip 

Let $\mathbf{P}({\R }^n, {\R }^n)$ be the space of ${\R }^n$ valued 
polynomials on ${\R }^n$. Then $\mathbf{P}({\R }^n, {\R }^n)$ $ = 
\mathrm{P}({\R }^n, \R) \otimes {\R }^n$, where $\mathrm{P}({\R }^n, \R )$ 
is the space of real valued polynomials on ${\R }^n$. We have a linear 
action 
\begin{displaymath}
\Theta : G \times \mathbf{P}({\R }^n, {\R }^n) \rightarrow 
\mathbf{P}({\R }^n, {\R }^n): (g, V) \mapsto ({\Phi }_g)_{\ast}V = 
{\Phi }_g (V \comp {\Phi }_{g^{-1}}). 
\end{displaymath}
When $V = p \otimes v$ for $p \in \mathrm{P}({\R }^n, \R)$ and $v \in {\R }^n$, 
\begin{displaymath}
{\Theta }_g(p \otimes v)(x) = ({\Phi}^{\dagger}_gp \otimes {\Phi }_{g^{-1}}(v))(x) = 
p({\Phi }_{g^{-1}}(x)) \otimes {\Phi }_{g^{-1}}v, 
\end{displaymath}
for every $(g,x) \in G\times {\R }^n$. To see this we compute. For every $h \in G$ 
\begin{align*}
{\Theta }_{gh}(p \otimes v)(x) & = p({\Phi }_{(gh)^{-1}}(x) \otimes 
{\Phi }_{(gh)^{-1}}v) 
= p({\Phi }_{h^{-1}}({\Phi }_{g^{-1}}(x))) \otimes {\Phi }_{h^{-1}}({\Phi }_{g^{-1}}v) \\
& = {\Phi }^{\dagger}_hp ({\Phi }_{g^{-1}}(x) \otimes 
{\Phi }_{h^{-1}}({\Phi }_{g^{-1}}v ) )= 
{\Theta }_h( ({\Phi }^{\dagger}_gp \otimes {\Phi }_{g^{-1}}v)(x) ) \\
& = {\Theta }_g ( {\Theta }_h(p \otimes v))(x).
\end{align*}
Let ${\mathrm{P}}_m({\R }^n, \R )$ be the space of homogeneous polynomials 
of degree $m$ on ${\R }^n$. Then $\mathrm{P}({\R }^n , \R ) = 
\sum_{m\ge 0} \oplus {\mathrm{P}}_m({\R }^n, \R )$. When $m =1$, one has 
${\mathrm{P}}_1({\R }^n, \R ) = ({\R }^n)^{\ast }$, the space of linear functions 
on ${\R }^n$. So ${\mathrm{P}}_1(({\R }^n)^{\ast }, \R ) = {\R }^n$. Consequently, 
\begin{align*}
{\mathbf{P}}_m({\R }^n , {\R }^n) & = 
{\mathrm{P}}_m({\R }^n , \R) \otimes {\R }^n = 
{\mathrm{P}}_m({\R }^n, \R) \otimes {\mathrm{P}}_1(({\R }^n)^{\ast }, \R ) \\
& \subseteq \sum^{\infty}_{k, \ell =0} {\mathrm{P}}_k({\R }^n, \R) 
\otimes {\mathrm{P}}_{\ell }(({\R }^n)^{\ast }, \R ) = \hspace{-5pt}
\sum^{\infty}_{\stackrel{k+\ell = 0}{k \ge 0 \, \, \& \, \, \ell \ge 0}} 
\hspace{-8pt}{\mathrm{P}}_{k+\ell }({\R }^n \oplus ({\R }^n)^{\ast }, \R ) \\
& = \mathrm{P}({\R }^n \oplus ({\R }^n)^{\ast }, \R ). 
\end{align*}
The above inclusion was used in \cite{cushman-sanders-white}. \medskip 

Let 
\begin{equation}
\Psi : G \times \mathrm{P}({\R }^n \oplus ({\R }^n)^{\ast }, \R ) 
\rightarrow \mathrm{P}({\R }^n \oplus ({\R }^n)^{\ast }, \R ) : 
(g, q) \mapsto {\Psi }_gq, 
\label{eq-s1two}
\end{equation}
where $({\Psi }_gq)(x, \xi ) = q({\Phi }_{g^{-1}}(x), {\Phi }^T_g\xi )$ for every  
$(x, \xi ) \in {\R }^n \times ({\R }^n)^{\ast }$. The 
mapping $\Psi $ is well defined, because ${\Phi }_{g^{-1}}: {\R }^n \rightarrow 
{\R}^n$ and ${\Phi }^T_g: ({\R }^n)^{\ast } \rightarrow ({\R }^n)^{\ast }$ are 
linear mappings. This implies that the mapping 
$(x, \xi) \mapsto q({\Phi }_{g^{-1}}(x), {\Phi }^T_g\xi )$ 
is a polynomial on ${\R }^n \oplus ({\R }^n)^{\ast }$, since 
$q$ is a polynomial on ${\R }^n \oplus ({\R }^n)^{\ast }$. For every 
$g \in G$ the map $q \mapsto {\Psi }_gq$, which sends  
$\mathrm{P}({\R }^n \oplus ({\R }^n)^{\ast }, \R )$ into itself, is linear. $\Psi $ is an action, because for every $g$, $h \in G$ and 
every $q \in \mathrm{P}({\R }^n \oplus ({\R }^n)^{\ast }, \R ) $ 
\begin{align*}
({\Psi }_{gh}q)(x, \xi ) & = q({\Phi }_{(gh)^{-1}}(x), {\Phi }^T_{gh}\xi ) 
= q ({\Phi }_{h^{-1}}({\Phi }_{g^{-1}}(x)), {\Phi }^T_h({\Phi }^T_g\xi )) \\
& = {\Psi }_hq ({\Phi }_{g^{-1}}(x), {\Phi }^T_g\xi ) =  
{\Psi }_g( {\Psi }_hq) (x, \xi ). 
\end{align*}

\noindent \textbf{Proposition 4} The subalgebra ${\mathrm{P}}({\R }^n \oplus 
({\R }^n)^{\ast } )^G$ of polynomials on ${\R }^n \oplus ({\R }^n)^{\ast }$, which are 
invariant under the $G$ action $\Psi $, is finitely generated. \medskip 

\noindent \textbf{Proof.} This follows from fact 2 because $G$ is a compact 
subgroup of $\Gl (n, \R )$ and the action $\Psi $ is linear. \hfill $\square $ \medskip 

Let ${\mathrm{P}}^{\mathrm{lin}}_{m+1}({\R }^n \oplus ({\R }^n)^{\ast })$ be the vector subspace of $\mathrm{P}({\R }^n \oplus ({\R }^n)^{\ast }, \R )$ consisting of homogeneous polynomials $q_{m+1}$ on ${\R }^n \oplus ({\R }^n)^{\ast }$ of degree $m+1$, which are linear in $({\R }^n)^{\ast }$. In other words, elements 
of $\mathrm{P}^{\mathrm{lin}}_{m+1}({\R }^n \oplus ({\R }^n)^{\ast }, \R )$ are of the form 
\begin{equation}
q_{m+1}(x,\xi ) = \sum_{i =1} q_{m,i}(x){\xi }_i = \sum_{|j|=1} q_{m,j}(x){\xi }^j, 
\label{eq-s1three}
\end{equation}  
where $q_{m,j} \in {\mathrm{P}}_m( {\R }^n, \R )$. Here $j = (j_i, \ldots , j_n) \in ({\Z }_{\ge 0})^n$ and $|j| = j_1 + \cdots + j_n$. 
Also $x \in {\R }^n$ and ${\xi }^j = {\xi }^{j_1}_1 \cdots {\xi }^{j_n}_n$ for 
$({\xi }_1, \ldots , {\xi }_n) \in ({\R }^n)^{\ast }$. For every $g \in G$   
\begin{align*}
({\Psi }_g q_{m+1})(x, \xi ) & = q_{m+1}({\Phi }_{g^{-1}}(x), {\Phi }^T_g\xi )  = \sum^n_{i=1}q_{m,i}({\Phi }_{g^{-1}}(x))({\Phi }^T_g\xi )_i \\
& = \sum_{|j| =1} q_{m,j}({\Phi }_{g^{-1}}(x))({\Phi }^T_g \xi )^j, 
\end{align*}
since $({\Phi }^T_g\xi )^j = ({\Phi }^T_g\xi )_j$ when $|j| =1$. So 
${\Psi }_gq_{m+1} \in \mathrm{P}^{\mathrm{lin}}_{m+1}
({\R }^n \oplus ({\R }^n)^{\ast }, \R )$. Thus $\mathrm{P}^{\mathrm{lin}}_{m+1}
({\R }^n \oplus ({\R }^n)^{\ast }, \R )$ is invariant under the action $\Psi $. 
Since $q_{m+1}(x,\xi )$ corresponds to the homogeneous polynomial vector field 
$\sum^n_{i=1} q_{m,i} \otimes {\xi }_i$ of degree $m$, we have proved \medskip 

\noindent \textbf{Corollary 4A} The vector space 
${\mathrm{P}}^{\mathrm{lin}}_{m+1}({\R }^n \oplus ({\R }^n)^{\ast })^G$ of $G$ invariant elements of ${\mathrm{P}}^{\mathrm{lin}}_{m+1}({\R }^n \oplus 
({\R }^n)^{\ast })$ is equal to the set $\mathcal{P}{\mathfrak{X}}_m({\R }^n)^G$ of homogeneous $G$ invariant polynomial vector fields on degree $m$ on ${\R }^n$. \medskip 

\noindent \textbf{Corollary 4B} The ${\mathrm{P}}({\R }^n, \R)^G$ module 
$\mathcal{P}\mathfrak{X}({\R }^n)^G$ of $G$ invariant polynomial 
vector fields on ${\R }^n$ is finitely generated. \medskip 

\noindent \textbf{Proof.} This follows from proposition 4 and corollary 
4A. \hfill $\square $ \medskip 

\noindent \textbf{Proof of theorem 1.} Since $G$ is compact and the 
action $\Psi $ (\ref{eq-s1two}) is linear, by Schwarz' theorem \cite{schwarz} 
the polynomial generators of 
${\mathrm{P}}^{\mathrm{lin}}({\R }^n \oplus ({\R }^n)^{\ast })^G$ generate the algebra of smooth $G$ invariant functions on ${\R }^n \oplus ({\R }^n)^{\ast }$, which are linear in $({\R }^n)^{\ast }$. Since element of 
${\mathrm{P}}^{\mathrm{lin}}({\R }^n \oplus ({\R }^n)^{\ast })^G$ correspond to 
$G$ invariant polynomial vector fields, the vector fields corresponding 
to these generators are a basis of the $C^{\infty}({\R }^n)^G$ module 
$\mathfrak{X}({\R }^n)^G$ of smooth $G$ invariant vector fields on ${\R }^n$. Thus the $C^{\infty}({\R }^n)^G$ module 
$\mathfrak{X}({\R }^n)^G$ is finitely generated. \hfill $\square $ \medskip 

Let 
\begin{equation}
\sigma : {\R }^n \rightarrow {\R }^k: x \mapsto \big( p_1(x), \ldots , p_k(x) \big) 
\label{eq-s1four}
\end{equation}
be the Hilbert map associated to the generators $p_1, \ldots , p_k$ of 
$\mathrm{P}({\R }^n,\R )^G$. Then $\Sigma = \sigma ({\R }^n)$ is a semialgebraic 
variety in ${\R }^k$. Moreover, $\sigma $ is the $G$ orbit mapping 
\begin{equation}
\pi : {\R }^n \rightarrow \Sigma = {\R }^n/G \subseteq {\R }^k: x \mapsto G\cdot x. 
\label{eq-s1five}
\end{equation}
The differential spaces $({\R }^n/G, C^{\infty}({\R }^n/G))$ and 
$(\Sigma , C^{\infty}({\R }^k)_{|\Sigma })$ are diffeomorphic. \medskip 

Let $\mathcal{P}\mathfrak{X}({\R }^n/G)$ be the set of all polynomial 
vector fields on the orbit space ${\R }^n/G$. In other words, 
$Y \in \mathcal{P}\mathfrak{X}({\R }^n/G)$ if and only if $Y$ is 
a derivation of the space of polynomials on ${\R }/G$, whose 
maximal integral curve through $y \in {\R }^n/G$ is an open interval, 
see \cite[prop. 3.2.6]{sniatycki13}. \medskip 

\noindent \textbf{Lemma 5} Let $X \in \mathcal{P}\mathfrak{X}({\R }^n)^G$. 
Then the vector field $Y = {\sigma }_{\ast } \comp X \comp {\sigma }^{\ast }$ lies 
in $\mathcal{P}\mathfrak{X}({\R }^n/G)$. \medskip 

\noindent \textbf{Proof.} Let $\overline{p}$ be a polynomial on ${\R }^n/G$, 
that is, there is a polynomial $P$ on ${\R }^k$ such that $P_{|{\R}^n/G} = p$. 
Then ${\sigma }^{\ast }\overline{p}$ is a $G$ invariant polynomial $q$ on 
${\R }^n$, for 
\begin{displaymath}
q(x) = {\sigma }^{\ast }\overline{p}(x) = \overline{p}(\sigma (x)) = 
(P \comp \sigma )(x). 
\end{displaymath}
For every $y \in {\R }^n/G$ one has $(Y\overline{p})(y) = 
{\sigma }_{\ast }(X(q))(y)$, since ${\sigma }^{\ast }\overline{p} = q$. Consequently, $Y$ maps a polynomial on ${\R }^n/G$ to a polynomial on ${\R }^n/G$. By definition $Y$ is a derivation, which is $\sigma $ related to the $G$ invariant polynomial vector field $X$ on ${\R }^n$. Because 
$({\R }^n/G, C^{\infty}({\R }^n/G))$ is a locally closed subcartesian differential 
space, $Y$ is a vector field on ${\R }^n/G$. \hfill $\square $\medskip 

Let $H$ be a compact subgroup of $G$. The orbit type ${\R }^n_{(H)}$ of the linear action $\Phi $ (\ref{eq-s1one}) is a locally trivial smooth fibration 
\begin{displaymath}
{\pi }_{|(H)}: {\R }^n_{(H)} \rightarrow \pi ({\R }^n_{(H)}) = 
{\overline{\R }}^n_{(H)} 
\end{displaymath}
with $G$ orbits as fibers, see lemma 2 of \cite{bates-cushman-sniatycki}. Here the $G$ orbit mapping $\pi $ (\ref{eq-s1five}) is 
the Hilbert mapping (\ref{eq-s1four}). Since the orbit space ${\R }^n/G$ is 
connected, there is a unique principal orbit type ${\R }^n_{(K)}$, which is 
open and dense. Moreover, ${\overline{\R}}^n_{(K)}$ is connected, see Duistermaat and Kolk \cite{duistermaat-kolk}. \medskip 

\noindent \textbf{Theorem 6} Every polynomial vector field on ${\R }^n/G$ is 
$\pi $ related to a $G$ invariant polynomial vector field on ${\R }^n$. \medskip 

\noindent \textbf{Proof.} Let $Y$ be a polynomial vector field on the orbit 
space ${\R }^n/G$. Then $Y_{|{\overline{\R}}^n_{(K)}}$ is a polynomial 
vector field on ${\overline{R}}^n_{(K)}$. We now take a closer look at the 
proof of theorem 1 of \cite{bates-cushman-sniatycki}. Let $U$ be an open subset of ${\R }^n/G$ for which 
a trivialization ${\tau }_{|{\pi }^{-1}(U)}: {\pi }^{-1}(U) \rightarrow U \times G$ of 
the locally trivial fibration ${\pi }_{|{\R}^n_{(K)}}: {\R }^n_{(K)} \rightarrow 
{\overline{\R }}^n_{(K)}$ is defined. There is a 
$G$ invariant vector field $\widetilde{Y}$ on $U \times G$ which is 
${\pi }_1: U \times G \rightarrow U:(u,g) \mapsto u$ related to $Y_{|U}$. Pulling 
$\widetilde{Y}$ back by the trivialization ${\tau }_{|{\pi }^{-1}(U)}$ yields 
a $G$ invariant vector field ${\widetilde{X}}_{{\pi}^{-1}(U)}$ on 
${\pi }^{-1}(U)$. Restricting ${\widetilde{X}}_{{\pi}^{-1}(U)}$ to 
${\tau }^{-1}_{{\pi }^{-1}(U)}$ to $U \times \{ e\} $  gives  
\begin{displaymath}
{\widetilde{X}}_{{\pi }^{-1}(U)}({\tau }^{-1}_{{\pi }^{-1}(U)}(U \times \{ e \})) = 
({\tau }^{-1}_{{\pi }^{-1}(U)})^{\ast }({\widetilde{Y}}_{|{U \times \{ e \}}}), 
\end{displaymath}
that is, 
\begin{equation}
{X}^{\vvee} = ({\tau }_{{\pi }^{-1}(U)})_{\ast } \big( {\widetilde{X}}_{{\pi }^{-1}(U)}
({\tau }^{-1}_{{\pi }^{-1}(U)}(U \times \{ e \}) \big) = 
{\widetilde{Y}}_{|(U \times \{ e \})}. 
\label{eq-s1six}
\end{equation} 
Pushing both sides of equation (\ref{eq-s1six}) forward by the projection 
mapping ${\pi }_1$ and using ${\pi }_1 \comp {\tau }_{{\pi }^{-1}(U)} = 
\pi $ gives ${\pi }_{\ast }{X}^{\vvee} = Y_{|U}$. Since ${\pi }_{|{\pi}^{-1}(U)}$ is a polynomial mapping and the vector field 
$Y_{|U}$ is polynomial, using equation (\ref{eq-s1six}) it follows that 
the vector field ${X}^{\vvee}$ is polynomial. Let $X_{{\pi }^{-1}(U)}(g \cdot x) = 
({\Phi }_g)_{\ast }X^{\vvee}(x)$ for every $g \in G$ and every 
$x \in {\tau }^{-1}_{{\pi }^{-1}(U)}(U \times \{ e \} )$. Then $X_{{\pi }^{-1}(U)}$ is 
a vector field on ${\pi }^{-1}(U)$, because 
\begin{displaymath}
\bigcup_{g \in G} {\Phi }_g({\tau }^{-1}_{{\pi }^{-1}(U)}(U \times \{ e \} )) = 
\bigcup_{g \in G}{\tau }^{-1}_{{\pi }^{-1}(U)}{\phi }_g(U \times \{ e \} ) = 
{\tau }^{-1}_{{\pi }^{-1}(U)}(U \times G) = {\pi }^{-1}(U),  
\end{displaymath} 
where $\phi : G \times (U \times G)  \rightarrow U\times G: (g,(u, h)) \mapsto 
(u,gh)$. Moreover, the vector field $X_{{\pi }^{-1}(U)}$ is polynomial, since 
${\Phi }_g$ is a linear mapping on ${\R }^n$ into itself for every $g \in G$ 
and $X^{\vvee}$ is a polynomial vector field. By construction the vector field $X_{{\pi }^{-1}(U)}$ is $\pi $ related to 
the vector field $Y_{|U}$. Since ${\pi }^{-1}(U)$ is an open subset of 
${\R }^n_{(K)}$, which in turn is an open subset of ${\R }^n$, the polynomial 
vector field $X_{{\pi }^{-1}(U)}$ extends to a unique $G$ invariant 
polynomial vector field $X$ on ${\R }^n$. The polynomial vector field 
$Y_{| {\overline{ {\R } }}^n_{(K)} }$ on ${\overline{\R}}^n_{(K)} \subseteq {\R }^k$ is 
$\pi $ related to the vector field $X_{{\pi }^{-1}(U)}$ on ${\pi }^{-1}(U)$, 
which is the domain of a trivialization of the bundle ${\pi }_{|{\R }^n_{(K)}}: 
{\R }^n_{(K)} \rightarrow {\overline{\R }}^n_{(K)}$. This implies that the 
vector fields $X$ and $Y$ are $\pi $ related, because 
${\R }^n_{(K)}$ is dense in ${\R }^n$. \hfill $\square $ \medskip 

\noindent \textbf{\large Ackowledgment}\medskip 

I would like to thank Prof. David Chillingworth for bringing references \cite{golubitsky-stewart-schaeffer} and \cite{po\`{e}naru} to my attention.

\end{document}